\documentclass{ajour}
\usepackage{amsmath,amssymb}

\newcommand{\cat}{\operatorname{cat}}
\renewcommand{\tfrac}{\textstyle\frac}
\newcommand{\R}{{\mathbb R}} 
\newcommand{\Rn}{{{\mathbb R}^n}} 
\newcommand{\C}{{\mathbb{C}}}

\newcommand{\irn}{\int\limits_{\Rn}}

\newcommand{\re}{\operatorname{Re}}

\renewcommand{\a }{\alpha }
\renewcommand{\b }{\beta }

\newcommand{\e }{\varepsilon }

\newcommand{\n }{\nabla }

\newcommand{\de}{\partial}
\newenvironment{Remark}{\begin{remark} \rm}{\end{remark}}

\begin{document}

\authorrunninghead{Cingolani, Secchi}
\titlerunninghead{Schr\"odinger equations with electromagnetic fields}

\title{\bf Semiclassical limit for
nonlinear Schr\"odinger equations with
electromagnetic fields}


\author{Silvia CINGOLANI \thanks{Supported by MURST}}
\affil{Dip. Inter. Matematica, Politecnico di Bari, Via E.Orabona 4,
70125 Bari, Italy}
\email{cingolan@dm.uniba.it}

\and

\author{Simone SECCHI \thanks{Supported by MURST, national project
\emph{Variational methods and nonlinear differential equations}.}}
\affil{SISSA, via Beirut 2/4, 34014 Trieste, Italy}
\email{secchi@ma.sissa.it}


\begin{article}


%
%
%
%
%
%

\section{Introduction}\label{sec:int}

Let us consider the  nonlinear Schr\"{o}dinger equation
\begin{equation}\label{eq:1.1}
ih \frac{\de \psi}{\de t} = \left( \frac{h}{i} \nabla -A(x)\right)^2
\psi +U(x)\psi - f(x,\psi), \quad x\in\Rn
\end{equation}
where $t \in \R$, $x \in \Rn$ $(n \geq 2)$.
The function $\psi(x,t)$ takes on complex values,
$h$ is the Planck constant, $i$ is the imaginary unit.
Here  $A\colon \Rn \to \Rn$ denotes a magnetic potential and
the Schr\"odinger operator is defined by
$$
\left( \frac{h}{i} \nabla -A(x)\right)^2 \psi
:= -h^2 \varDelta \psi - \frac{2h}{i} A \cdot \nabla \psi + |A|^2 \psi -
\frac{h}{i} \psi\, \operatorname{div}A\,  .
$$
Actually, in general dimension $n \geq 2$,
the magnetic field $B$ is a $2$-form where
$B_{i,j}= \partial_j A_k - \partial_k A_j$ ;
in the case $n=3$, $B= \operatorname{curl}A$.
The function $U\colon \Rn \to \R$ represents an electric potential.
In the sequel, for the sake of simplicity, we limit ourselves
to the particular case in which $f(x,t)=K(x)|t|^{p-1}t$, with $p>1$ if
$n=2$ and $1<p<\frac{n+2}{n-2}$ if $n\geq 3$.

It is now well known that the
nonlinear Schr\"{o}dinger equation \eqref{eq:1.1}
arises from a perturbation approximation for strongly nonlinear
dispersive wave systems. Many papers are devoted to the nonlinear
Schr\"{o}dinger equation and its solitary wave solutions.

In this paper we seek for standing wave solutions to
\eqref{eq:1.1}, namely waves of the form
$\psi (x,t)=e^{-iEth^{-1}}u(x)$ for some function
$u\colon \Rn \longrightarrow \C$.
Substituting this
{\it ansatz} into \eqref{eq:1.1}, and denoting for
convenience $\e=h$, one is led to solve the complex equation in $\Rn$
\begin{equation}
\left( \tfrac{\e}{i} \nabla -A(x)\right)^2 u +
(U(x)-E) u = K(x)|u|^{p-1}u \,.
\tag{NLS}
\end{equation}
Renaming $V(x)+1=U(x)-E$, we assume from now on that $1+V$ is
strictly positive on the whole $\Rn$.
Moreover, by an obvious change of variables, 
the problem becomes that of finding some function
$u\colon\Rn\to\C$ such that
\begin{equation}\label{S}
\left( \tfrac{\nabla}{i} -A(\e x)\right)^2 u + u +
V(\e x)u = K(\e x)|u|^{p-1}u, \quad x\in\Rn \,.
\tag{$S_\e$}
\end{equation}
Concerning nonlinear Schr\"odinger equation with external magnetic field,
we firstly quote a paper by Esteban and Lions \cite{EL}, where
concentrations and compactness arguments are applied to solve some
minimization problems associated to $(S_\e)$ under
suitable assumptions on the magnetic field.

The purpose of this paper is 
to study the time--independent nonlinear
Schr\"{o}dinger equation $(S_\e)$ in the semiclassical limit.
This seems a very interesting problem since
the Correspondence's Principle establishes that
Classical Mechanics is, roughly speaking, contained in
Quantum Mechanics.
The mathematical transition is obtained letting to zero
the Planck constant ($\e \to 0$)
and solutions $u(x)$ of $(S_\e)$
which exist for small value of $\e$ are
usually referred as semi-classical ones (see \cite{RS}).

We remark that in the
linear case, Helffer {\it et al.} in
\cite{hel1,helS} have studied the asymptotic
behavior of the eigenfunctions
of the  Schr\"{o}dinger operators with magnetic fields
in the semiclassical limit.
Note that in these papers
the {\sl wells} of the  Schr\"{o}dinger operators with magnetic fields
are the same as those without magnetic field,
so that one doesn' t `see'  the magnetic field in the
definition of the {\sl well}. See also \cite{brum}
for generalization of the results by \cite{helS}
for potentials which degenerate at infinity.

In the case $A=0$, (no magnetic field),
a recent extensive literature is devoted to
study the time--independent nonlinear
Schr\"{o}dinger equation $(S_\e)$ in the semi-classical limit.
We shortly recall the main results in literature.
The first paper is due to Floer and Weinstein which
investigated the one-dimensional nonlinear Schr\"{o}dinger equation
(with $K(x)=1$)  and gave a description of the limit behavior of $u(x)$
as $\e \to 0$.
Really they proved that if the potential $V$ has a non-degenerate
critical point, then $u(x)$ concentrates near this
critical point, as $\e \to 0$.

Later, other authors proved that this problem is really local in nature and
the presence of an isolated critical point of the potential $V$
(in the case $K(x)=1$)
produces a semi-classical solution $u(x)$ of $(S_\e)$ which
concentrates near this point.
Different approaches are used to cover different cases
(see \cite{ABC,dpf,Li,oh,xw}).
Moreover when $V$ oscillates, the existence of multibumps solutions has also
been studied in \cite{amberti,CN,dpf1,gui}.
Furthermore  multiplicity results are obtained in 
 \cite{CL1,CL2} for potentials $V$ having a set of degenerate
global minima and recently in \cite{ams},
for potentials $V$ having a set of
critical points, not necessarily global minima.

\smallskip
A natural answer arises:
{\sl how does the presence of an external magnetic field 
influence the existence and the
concentration behavior of
standing wave solutions to \eqref{eq:1.1} in the semi-classical limit}?
\smallskip

A first result in this direction is contained in \cite{ku}
where Kurata has proved the existence of
least energy solutions to $(S_\e)$ for any $\e >0$,
under some assumptions linking the magnetic field
$B= (B_{i,j})$ and the electric potential $V(x)$.
The author also investigated the semi-classical limit of the found
least energy solutions and
showed a concentration phenomenon near global
minima of the electric potential in the case $K(x) = 1$
and  $|A|$ is small enough.

Recently in \cite{cingolani}, Cingolani
obtained a multiplicity result of semi-classical standing waves
solutions to $(S_\e)$,
relating the number of solutions to $(S_\e)$ to
the {\sl richness} of a set $M$ of global minima 
of an auxiliary function $\Lambda$ defined by setting
\[
\Lambda (x)=\frac{(1+V(x))^\theta}{K(x)^{-2/(p-1)}},\qquad \theta =
\frac{p+1}{p-1}-\frac{n}{2},
\]
(see $(\ref{def:auxiliary})$ in section 4 for details) depending on $V(x)$ and $K(x)$.
We remark that, if $K(x)=1$ for any $x \in \Rn$, global
minima of $\Lambda$ coincides with global minima of $V$.
The variational approach, used in \cite{cingolani},
allows to deal with unbounded potential $V$
and does not require assumptions on the magnetic field.
However this approach works only near global minima of
$\Lambda$.

In the present paper we deal with the more
general case in which the auxiliary function
$\Lambda$ has a manifold $M$ of stationary points,
not necessarily global minima. 
For bounded magnetic potentials $A$, 
we are able to prove a multiplicity result of semi-classical
standing waves of $(S_\e)$,
following the new perturbation approach
contained in the recent paper \cite{ams} due to
Ambrosetti, Malchiodi, Secchi.

Now we briefly describe the proof of the result.
First of all, we highlight that solutions of $(S_\e)$ naturally
appear as {\em orbits}: in fact, equation $(S_\e)$ is invariant under the
multiplicative action of $S^1$. Since there is no danger of confusion, we
simply speak about solutions.
The complex--valued solutions to $(S_\e)$ are
found near least energy solutions of the equation \begin{equation}\label{risc}
\bigg(\frac{\nabla}{i} - A(\e \xi)\bigg)^2u + u + V(\e \xi)u=K(\e\xi)|u|^{p-1}u.
\end{equation}
where $\e \xi$ is in a neighborhood of $M$.
The least energy of (\ref{risc}) have the form
\begin{equation}\label{def:zed}
z^{\e\xi ,\sigma} \colon x\in\Rn \mapsto e^{i\sigma +
i A(\e\xi)\cdot x}
\left(\frac{1+V(\e\xi)}{K(\e\xi)}\right)^\frac{1}{p-1} 
U((1+V(\e\xi))^{1/2} (x-\xi))
\end{equation}
where $\e \xi$ belongs to $M$ and $\sigma \in [0, 2 \pi]$.
As in \cite{ams}, the proof relies on a suitable
finite dimensional reduction,
and critical points of the Euler functional $f_\e$ associated to problem $(S_\e)$
are found near critical point of a finite dimensional
functional $\Phi_\e$ which is defined on a suitable neighborhood of $M$.
This allows to use Lusternik-Schnirelman category in the case $M$ is a
set of local maxima or minima of $\Lambda$.
We remark that the case of maxima cannot be handled by using
direct variational arguments as in \cite{cingolani}.

Moreover in the case $M$ is  a set of critical points non-degenerate
in the sense of Bott (see \cite{bott}) we are able to prove the existence of (at least) cup long of 
$M$ solutions
concentrating near points of $M$. For the definition of the cup long, refer to section 5.

Firstly we present a special case of our results.

\begin{theorem}\label{th:nond}
Assume that
\begin{itemize}
\item[$(K1)$] \, $K\in L^\infty (\Rn)\cap C^2(\Rn)$ is strictly positive and $K''$ is
bounded; 
\item[$(V1)$] \, $V\in L^\infty (\Rn)\cap C^2(\Rn)$ satisfies $\inf_{x\in\Rn}
(1+V(x)) > 0$, and $V''$ is bounded;
\item[$(A1)$] \, $A\in L^\infty (\Rn,\Rn)\cap C^1(\Rn,\Rn)$, and
the jacobian $J_A$ of $A$ is globally bounded in $\Rn$.
\end{itemize}

\noindent
If the auxiliary function $\Lambda$ has
a non-degenerate critical point $x_0\in \Rn$, then
for $\e > 0$ small, the problem $(S_\e)$ has at least a
(orbit of) solution concentrating near $x_0$.
\end{theorem}

Actually, we are able to prove the following generalization.

\begin{theorem}\label{th:main}
As in Theorem 1.1, assume again (K1), (V1) and (A1) hold.

\noindent
If the auxiliary function $\Lambda$ has
a smooth, compact, non-degenerate manifold of critical points $M$,
then for $\e > 0$ small, the problem $(S_\e)$ has at least
$\ell (M)$ (orbits of) solutions concentrating near points of $M$.
\end{theorem}

Finally we point out that the presence of
an external magnetic field
produces a phase in the complex wave which depends on
the value of $A$ near $M$.
Conversely the presence of $A$ does not seem to influence the location of
the peaks of the modulus of the complex wave. Although we will not deal with this problem, 
we believe that in order to have a 
local $C^2$ convergence of the solutions, some assumption about the smallness of the 
magnetic potential $A$ should be added, as done in \cite{ku} for minima of $V$.

\bigskip
Finally we point out that
Theorem~\ref{th:nond} and Theorem~\ref{th:main} hold for problems involving
more general nonlinearities. See Remark \ref{nonl} in the last section.

\vspace{10pt}

\begin{center}\textbf{Notation}\end{center} 
\begin{enumerate}
\item The complex conjugate of any number $z\in\C$ will be denoted by $\bar z$.
\item The real part of a number $z\in\C$ will be denoted by $\re z$.
\item The ordinary inner product between two vectors $a,b\in\Rn$ will be denoted by $a \cdot 
b$.
\item From time to time, when no confusion can arise, we omit the symbol $dx$ in integrals over
$\Rn$.
\item $C$ denotes a generic positive constant, which may vary inside a chain of inequalities.
\item
We use the Landau symbols. For example $O(\e)$ is a generic function such that 
$\limsup\limits_{\e\to 0} \frac{O(\e)}{\e} <
\infty$, and $o(\e)$ is a function such that $\lim\limits_{\e\to 0} \frac{o(\e)}{\e}=0$.
\end{enumerate}

\section{The variational framework}
We work in the real Hilbert space $E$ obtained as the completion of
$C_0^\infty (\Rn,\C)$ with respect to the norm associated to the inner product
\[
\langle u\mid v\rangle = \operatorname{Re} \int_\Rn \nabla u \cdot \overline{\nabla v} + u \bar v 
..
\]
Solutions to ($S_\e$) are, under some conditions we are going to point out, critical points of the 
functional formally defined on $E$ as
\begin{multline}\label{feps}
f_\e (u)=\frac{1}{2}\int_\Rn \bigg(\bigg| \bigg(\frac{1}{i}\nabla - A(\e x)\bigg)u\bigg|^2 
+|u|^2+V(\e x)|u|^2\bigg)\, dx \\- \frac{1}{p+1}\int_\Rn K(\e x)|u|^{p+1}\, dx.
\end{multline}
In what follows, we shall assume that the functions $V$, $K$ and $A$
satisfy the following assumptions:
\begin{itemize}
\item[(K1)] $K\in L^\infty (\Rn)\cap C^2(\Rn)$ is strictly positive and $K''$ is
bounded; 
\item[(V1)] $V\in L^\infty (\Rn)\cap C^2(\Rn)$ satisfies $\inf_{x\in\Rn}
(1+V(x)) > 0$, and $V''$ is bounded;
\item[(A1)] $A\in L^\infty (\Rn,\Rn)\cap C^1(\Rn,\Rn)$, and
the jacobian $J_A$ of $A$ is globally bounded in $\Rn$.
\end{itemize}
Indeed,
\begin{multline*}
\int_\Rn \bigg(\bigg| \bigg(\frac{1}{i}\nabla - A(\e x)\bigg)u\bigg|^2 \bigg)\, dx
= \\ = \int_\Rn \bigg(|\nabla u|^2 + |A(\e x)u|^2 -2
\operatorname{Re} ( \textstyle\frac{\nabla u}{i} \cdot A(\e x)\overline{u}) \bigg)\, dx,
\end{multline*}
and the last integral is finite thanks to the Cauchy--Schwartz
inequality and the boundedness of $A$.

\noindent It follows that $f_\e$ is actually well-defined on $E$.

In order to find possibly multiple critical points of \eqref{feps},
we follow the approach of \cite{ams}. In our context, we need to find
complex--valued solutions, and so some further remarks are due.

Let $\xi\in\Rn$, which will be fixed suitably later on: we look for solutions to \eqref{S} ``close'' to 
a particular solution of the equation
\begin{equation}\label{eqlimite}
\bigg(\frac{\nabla}{i} - A(\e \xi)\bigg)^2u + u + V(\e\xi)u=K(\e\xi)|u|^{p-1}u.
\end{equation}
More precisely, we denote by $U_c \colon \Rn \to \C$ a least--energy solution
to the scalar problem
\begin{equation}\label{eq:4}
-\varDelta U_c + U_c + V(\e\xi)U_c = K(\e\xi) |U_c|^{p-1}U_c
\mbox{\quad in }\Rn.
\end{equation}
By energy comparison (see \cite{ku}), one has that
\[
U_c (x)=e^{i\sigma} U^\xi (x-y_0)
\]
for some choice of $\sigma \in [0,2\pi]$ and $y_0\in\Rn$, where
$U^\xi\colon\Rn\to\R$ is the unique solution of the problem
\begin{equation}
\begin{cases}
-\varDelta U^\xi + U^\xi + V(\e\xi) U^\xi = K(\e\xi)| U^\xi|^{p-1} U^\xi  \\
U^\xi (0) = \max_\Rn U^\xi\\
U^\xi > 0.
\end{cases}
\end{equation}
If $U$ denotes the unique solution of
\begin{equation}
\begin{cases}
-\varDelta U + U = U^p \mbox{\quad in }\Rn\\
U(0) = \max_\Rn U \\
U > 0,
\end{cases}
\end{equation}
then some elementary computations prove that
$U^\xi (x)=\alpha(\e\xi)U(\beta(\e\xi)x)$,
where
\begin{eqnarray*}
\alpha (\e\xi) &= \left(\frac{1+V(\e\xi)}{K(\e\xi)}\right)^\frac{1}{p-1} \\
\beta(\e\xi)&=(1+V(\e\xi))^{1/2}.
\end{eqnarray*}
It is easy to show, by direct computation, that the function
$u(x)= e^{i A(\e\xi)\cdot x} U_c (x)$ actually solves \eqref{eqlimite}.

For $\xi \in \Rn$ and $\sigma\in [0,2\pi]$, we set
\begin{equation}\label{def:zeta}
z^{\e\xi ,\sigma} \colon x\in\Rn \mapsto e^{i\sigma + i A(\e\xi)\cdot 
x}\alpha(\e\xi)U(\beta(\e\xi)(x-\xi)).
\end{equation}

Sometimes, for convenience, we shall identify $[0,2\pi]$ and $S^1\subset \C$, through $\eta = 
e^{i\sigma}$.

Introduce now the functional $F^{\e\xi,\sigma}\colon E \to \R$ defined by
\begin{multline*}
F^{\e\xi,\sigma}(u)=\frac{1}{2}\int_\Rn  \bigg(\bigg|\bigg(\frac{\nabla u}{i} -
A(\e\xi)u \bigg)\bigg|^2 + |u|^2 + V(\e\xi)|u|^2 \bigg) \, dx \\ - \frac{1}{p+1}
\int_\Rn K(\e\xi)|u|^{p+1}\, dx,
\end{multline*}
whose critical point correspond to solutions of \eqref{eqlimite}.

The set
\[
Z^\e =\{{z^{\e\xi,\sigma}}\mid \xi \in \Rn \, \land \,
\sigma \in [0,2\pi]\} \simeq S^1 \times \Rn
\]
is a regular manifolds of critical points for the functional $F^{\e\xi,\sigma}$.

It follows from elementary differential geometry that
\begin{multline*}
T_{z^{\e\xi,\eta}}Z^\e = \operatorname{span}_{\R}
\{\tfrac{\de}{\de \sigma}
{{{{{z^{\e\xi,\sigma}}}}}},\tfrac{\de}{\de \xi_1}
{{z^{\e\xi,\sigma}}},\dots,\tfrac{\de}{\de \xi_n}
{{{z^{\e\xi,\sigma}}}}\}=\\ =\operatorname{span}_{\R}
\{i{{{z^{\e\xi,\sigma}}}},\tfrac{\de}{\de \xi_1}
{{{z^{\e\xi,\sigma}}}},\dots,\tfrac{\de}{\de \xi_n}{{z^{\e\xi,\sigma}}}\},
\end{multline*}
where we mean by the symbol $\operatorname{span}_{\R}$ that all
the linear combinations must have real coefficients.

We remark that, for $j=1,\dots,n$,
\begin{multline*}
\frac{\de}{\de \xi_j}{{z^{\e\xi,\sigma}}}=-\frac{\de}{\de x_j}{{z^{\e\xi,\sigma}}} + O(\e|\nabla 
V(\e\xi)|)+\\
+i\alpha (\e\xi)e^{i A(\e\xi)\cdot x + i\sigma} U(\beta(\e\xi)(x-\xi))\left(\frac{\de}{\de 
\xi_j}\left(A(\e\xi)\cdot x\right)
+A_j(\e\xi)\right)=\\
=-\frac{\de}{\de x_j} {{z^{\e\xi,\sigma}}} + O(\e|\nabla V(\e\xi)|+\e|J_A 
(\e\xi)|)+i{{{z^{\e\xi,\sigma}}}}
A_j(\e\xi),
\end{multline*}
so that
\[
\frac{\de}{\de \xi_j}{{z^{\e\xi,\sigma}}} =
- \frac{\de}{\de x_j}{{z^{\e\xi,\sigma}}} +
i {z^{\e\xi,\sigma}}A_j(\e\xi)+ O(\e).
\]
Collecting these remarks,
we get that any
$\zeta \in T_{{{{{z^{\e\xi,\sigma}}}}}}Z^\e$ can be written as
$$
\zeta = i\ell_1 {{{z^{\e\xi,\sigma}}}}+ \sum_{j=2}^{n+1} \ell_j
\frac{\de}{\de x_{j-1}}{{{{z^{\e\xi,\sigma}}}}}+O(\e)$$
for some real coefficients $\ell_1,\ell_2,\dots,\ell_{n+1}$.

The next lemma shows that $\nabla f_\e ({{{z^{\e\xi,\sigma}}}})$ gets small when
$\e\to 0$.

\begin{lemma}\label{lem:1}
For all $\xi\in \Rn$, all $\eta \in S^1$ and all $\e>0$ small, one has that 
$$ 
\|\nabla f_\e ({z^{\e\xi,\sigma}})\|\leq C\left(\e |\nabla 
V(\e\xi)|+\e |\nabla K(\e\xi)|+\e |J_A (\e\xi)|+\e 
|\operatorname{div} A(\e\xi)|+\e^{2}\right),  $$ 
for some constant $C>0$.
 \end{lemma} 
\begin{proof} 
From 
\begin{multline}\label{eq:7}
f_{\e}(u)  =  F^{\e\xi,\eta}(u)+\frac{1}{2}\int_\Rn \bigg( \bigg| \frac{\nabla
u}{i} - A(\e x)u \bigg|^2 - \bigg| \frac{\nabla u}{i}-A(\e\xi)u
\bigg|^2 \biggr)  +\\
+ \frac{1}{2}\int_{\Rn}\left[V(\e x)-V(\e\xi)\right]u^{2}
-\frac{1}{p+1}\int_\Rn [K(\e x)-K(\e\xi)]|u|^{p+1}
\end{multline}
 and since $z^{\e\xi,\eta}$ is a
critical point of $F^{\e\xi,\eta}$, one has (with $z=z^{\e\xi,\eta}$)
\begin{multline*} 
\langle \nabla f_\e (z)\mid v \rangle = \re \int_\Rn \bigg(
\frac{1}{i}\nabla - A(\e\xi)\bigg)z \cdot (A(\e\xi)-A(\e x))\bar v 
\\+ \re\int_\Rn (A(\e\xi)-A(\e x))z \cdot \overline{\bigg( \frac{1}{i}\nabla -
A(\e\xi)\bigg)v}  + \\ \re\int_\Rn (A(\e\xi)-A(\e x))z \cdot
(A(\e\xi)-A(\e x))\bar v \\
+ \re\int_\Rn (V(\e x)-V(\e\xi))z\bar v  - \re\int_\Rn (K(\e
x)-K(\e\xi))|z|^{p-2}z\bar v \\ = \e \operatorname{Re} \int_\Rn \frac{1}{i}
(\operatorname{div}\, A(\e x)) z \bar v  + 2 \operatorname{Re} \int_\Rn
(A(\e\xi)-A(\e x))z \cdot \overline{\bigg( \frac{1}{i} \nabla -
A(\e\xi)\bigg)v} \\+\re\int_\Rn (V(\e x)-V(\e\xi))z\bar v  -
\re\int_\Rn (K(\e x)-K(\e\xi))|z|^{p-2}z\bar v 
\end{multline*}

\noindent
From the assumption that $|D^{2}V(x)|\leq {\rm const.}$ one infers 
\[ 
|V(\e x)-V(\e\xi)|\leq \e |\nabla V(\e\xi)|\cdot
|x-\xi|+c_{1}\e^{2}|x-\xi|^{2}.
\]
This implies 
\begin{eqnarray} 
 \int_{\Rn}|V(\e x)-V(\e\xi)|^{2}{{{z^{\e\xi,\sigma}}}}^{2}
 &\leq&  c_1\e^{2}|\nabla V(\e\xi)|^{2}\int_{\Rn}|x-\xi|^{2}z^{2}(x-\xi) + \nonumber \\ 
 & & c_{2}\e^{4}\int_{\Rn}|x-\xi|^{4}z^{2}(x-\xi). 
\label{eq:1.3} 
\end{eqnarray} 
A direct calculation yields 
\begin{eqnarray*} 
\int_{\Rn}|x-\xi|^{2}z^{2}(x-\xi) & = & 
\a^{2}(\e\xi)\int_{\Rn}|y|^{2}U^{2}(\b(\e\xi) y)dy \\ 
& = & \alpha (\e\xi)^{2} \beta (\e\xi)^{-n-2}\int_{\Rn}|y'|^{2}U^{2}(y')dy'\leq c_{3}. 
\end{eqnarray*} 
{From} this (and a similar calculation for the last integral in the 
above formula) one infers 
\begin{equation} 
\int_{\Rn}|V(\e x)-V(\e\xi)|^{2}|z^{\e\xi,\sigma}|^{2}
\leq c_{4}\e^{2}|\nabla
V(\e\xi)|^{2} + c_{5}\e^{4}. 
\label{eq:1.4} 
\end{equation} 
Of course, similar estimates hold for the terms involving $K$.
It then follows that
\[
\|\nabla f_\e (z^{\e\xi,\eta})\| \leq C(\e | \operatorname{div}A(\e\xi)| + \e
| \nabla V (\e\xi)| + \e | J_A (\e\xi)| + \e^2),
\]
and the lemma is proved.
  \end{proof} 

\section{The invertibility of $D^2f_\e$ on $(TZ^\e)^\bot$}

To apply the perturbative method, we need to exploit some non--degeneracy properties of
the solution ${z^{\e\xi,\sigma}}$ as a critical point of $F^{\e\xi,\sigma}$.

Let $L_{\e,\sigma,\xi}\colon (T_{{{z^{\e\xi,\sigma}}}}Z^\e)^\bot
\longrightarrow
(T_{{z^{\e\xi,\sigma}}}Z)^\bot$ be the operator defined by \[\langle
L_{\e,\sigma,\xi} v \mid w \rangle =
D^2 f_\e ({{{z^{\e\xi,\sigma}}}})(v,w)\] for all
$v,w \in (T_{{z^{\e\xi,\sigma}}}Z^\e)^\bot$.

The following elementary result will play a fundamental role in the present
section.

\begin{lemma}
Let $M \subset \Rn$ be a bounded set. 
Then there exists a constant $C>0$ such that for all $\xi\in M$ one has
\begin{equation}\label{norme:equiv}
\int_\Rn \bigg| \bigg( \frac{\nabla}{i} - A(\xi) \bigg) u \bigg|^2 + |u|^2 \geq C \int_\Rn (|\nabla 
u|^2 + |u|^2) \qquad\forall u \in E.
\end{equation}
\end{lemma}

\begin{proof}
To get a contradiction, we assume on the contrary the existence of a sequence
$\{ \xi_n\}$ in $M$ and a sequence $\{u_n\}$ in $E$ such that $\|u_n\|_E = 1$
for all $n\in\mathbb{N}$ and 
\begin{equation}\label{eq:contradiction}
\lim_{n\to + \infty} \bigg[ \int_\Rn \bigg| \bigg( \frac{\n}{i} - A(\xi)
\bigg)u_n \bigg|^2 + \int_\Rn |u_n|^2 \bigg] = 0.
\end{equation}
In particular, $u_n \to 0$ strongly in $L^2 (\Rn,\C)$. Moreover, since $M$
is bounded, we can assume also $\xi_n \to \xi^* \in \overline{M}$ as
$n\to\infty$. From 

\begin{multline*}
\int_\Rn \bigg| \bigg( \frac{\n}{i} - A(\xi_n)
\bigg)u_n \bigg|^2 = \\
 \int_\Rn \left( |\n u_n|^2 + |A(\xi_n)|^2 |u_n|^2 - 2 \re
\frac{1}{i} \n u_n \cdot A(\xi_n) \overline{u_n} \right)
\end{multline*}
we get 
\begin{eqnarray*}
\lim_{n \to + \infty} \int_\Rn | \n u_n |^2 &=& 1 \\
\lim_{n \to + \infty} \re \int_\Rn \tfrac{1}{i} \n u_n \cdot A(\xi_n) \overline{u_n} &=& \frac{1}{2} \,.
\end{eqnarray*}
Therefore,
\begin{eqnarray*}
\limsup_{n\to\infty} \int_\Rn |\n u_n| \, |A(\xi_n)| \, |u_n| &\geq&
\limsup_{n\to\infty} \left| \int_\Rn \textstyle\frac{1}{i} \n u_n \cdot
A(\xi_n) \overline{u_n} \right| \\
&\geq& \limsup_{n\to\infty} \re \int_\Rn \tfrac{1}{i} \n u_n \cdot A(\xi_n)
\overline{u_n} = \frac{1}{2} \,.
\end{eqnarray*}
From this we conclude that
\begin{eqnarray*}
\frac{1}{2} &\leq& |A(\xi^*)| \limsup_{n\to\infty} \| \n u_n \|_{L^2} \|u_n
\|_{L^2} \leq \\
&\leq& |A(\xi^*)| \limsup_{n\to\infty} \|u_n \|_{L^2} = 0,
\end{eqnarray*}
which is clearly absurd. This completes the proof of the lemma.
\end{proof}

At this point we shall prove the following result:

\begin{lemma}\label{lemma:3.1}
Given $\bar\xi >0$, there exists $C > 0$ such that for $\e$ small enough one
has
\begin{equation}\label{eq:posdef}
|\langle L_{\e,\sigma,\xi}v\mid v \rangle | \geq C \|v\|^2, \quad \forall |\xi|
\leq \bar\xi, \;\forall\sigma\in [0,2\pi],\;\forall v\in
(T_{{z^{\e\xi,\sigma}}}Z^\e)^\bot.
\end{equation}
\end{lemma}

\begin{proof}
We follow the arguments in \cite{ams}, with some minor modifications due to the presence of 
$A$. Recall that $$T_{{z^{\e\xi,\sigma}}}Z=\operatorname{span}_{\R} \{ \frac{\de}{\de
\xi_1}{z^{\e\xi,\sigma}},\dots ,\frac{\de}{\de
\xi_n}{z^{\e\xi,\sigma}},{i z^{\e\xi,\sigma}} \}.$$
Define 
$$
\mathcal{V}=\operatorname{span}_{\R} \{ \tfrac{\de}{\de
x_1}{z^{\e\xi,\sigma}},\dots ,\tfrac{\de}{\de
x_n}{z^{\e\xi,\sigma}},{z^{\e\xi,\sigma}} ,i{z^{\e\xi,\sigma}}\}.
$$
As in \cite{ams}, it suffices to prove \eqref{eq:posdef} for all
$v\in\operatorname{span}_\R \{{{z^{\e\xi,\sigma}}},\phi\}$,
where $\phi \perp \mathcal{V}$.
More precisely, we shall prove that for some constants $C_1 >0$, $C_2 > 0$,
for all $\e$ small enough and all $|\xi|\leq \bar\xi$ the following hold:
\begin{equation}\label{eq:13}
\langle L_{\e,\sigma,\xi}{{z^{\e\xi,\sigma}}}\mid {z^{\e\xi,\sigma}} \rangle \leq
-C_1 < 0,
\end{equation}
\begin{equation}\label{eq:14}
\langle L_{\e,\sigma,\xi}\phi\mid \phi \rangle
\geq C_2 \|\phi\|^2 \quad\forall\phi\perp\mathcal{V}.
\end{equation}
For the reader's convenience, we reproduce here the expression for the second
derivative of $F^{\e\xi,\sigma}$:
\begin{multline*}
D^2 F^{\e\xi,\sigma}(u)(v,v)=\int_\Rn \bigg| \bigg(
\frac{\nabla}{i}-A(\e\xi)\bigg)v\bigg|^2 + |v|^2 + V(\e\xi)|v|^2 \\
- K(\e\xi)\bigg[(p-1) \re \int_\Rn |u|^{p-3} \re (u \bar v) u \bar v  +
\int_\Rn |u|^{p-1}|v|^2 \bigg].
\end{multline*}
Moreover, since ${{z^{\e\xi,\sigma}}}$ is a solution of \eqref{eqlimite}, we
immediately get
\[
\int_\Rn
\bigg(  \bigg| \bigg( \frac{\nabla}{i}-A(\e\xi)\bigg){{{{z^{\e\xi,\sigma}}}}}\bigg|^2
+ V(\e\xi)|{z^{\e\xi,\sigma}}|^2 +|{{z^{\e\xi,\sigma}}}|^2
 \bigg)  = K(\e\xi)\int_\Rn
|{z^{\e\xi,\sigma}}|^{p+1}.
\]
From this it follows readily that we can find some $c_0 > 0$ such that for all
$\e > 0$ small, all $|\xi| \leq \bar \xi$ and all $\sigma\in [0,2\pi]$ it
results
\begin{equation} D^2
F^{\e\xi,\sigma}({{z^{\e\xi,\sigma}}}) ({z^{\e\xi,\sigma}},{z^{\e\xi,\sigma}})<c_0 <0.
\end{equation}
Recalling \eqref{eq:7}, we find
\begin{multline*}
\langle L_{\e,\sigma,\xi}{{z^{\e\xi,\sigma}}}\mid {z^{\e\xi,\sigma}} \rangle =D^2
F^{\e\xi,\sigma}({{z^{\e\xi,\sigma}}})({z^{\e\xi,\sigma}},{z^{\e\xi,\sigma}})+ \\ +\int_\Rn [V(\e 
x)-V(\e\xi)]|{{z^{\e\xi,\sigma}}}|^2
-\int_\Rn [K(\e x)-K(\e\xi)]|{{z^{\e\xi,\sigma}}}|^{p+1}
+\\ \int_\Rn
\bigg( \bigg| \bigg(\frac{\nabla}{i}-A(\e x)\bigg){{{z^{\e\xi,\sigma}}}}\bigg|^2 - \bigg|
\bigg(\frac{\nabla}{i}-A(\e \xi)\bigg){{{z^{\e\xi,\sigma}}}}\bigg|^2
\bigg) \,.
\end{multline*}
It follows that
\begin{multline}
\langle L_{\e,\sigma,\xi}{{z^{\e\xi,\sigma}}}\mid {z^{\e\xi,\sigma}} \rangle \leq D^2
F^{\e\xi,\sigma}({{z^{\e\xi,\sigma}}})({z^{\e\xi,\sigma}},{z^{\e\xi,\sigma}}) + \\ +c_1
\e |\nabla V(\e\xi)|+c_2\e |\nabla K (\e\xi)| + c_3 \e |J_A (\e\xi)| + c_4
\e^2. \end{multline}
Hence \eqref{eq:13} follows. The proof
of \eqref{eq:14} is more involved. We first prove the following
claim.

\vspace{5pt}

{\em Claim.} There results \begin{equation}\label{eq:17}
D^2 F^{\e\xi}({z^{\e\xi,\sigma}})(\phi,\phi) \geq c_1 \|\phi\|^2 \quad
\forall\phi\perp\mathcal{V}.
\end{equation}

\vspace{3pt}

Recall that the complex ground state $U_c$ introduced in \eqref{eq:4} is a
critical point of mountain-pass type for the corresponding energy functional
$J\colon E \longrightarrow \R$ defined by
\begin{equation}
J(u)=\frac{1}{2}\int_\Rn ( |\nabla u|^2 + |u|^2 + V(\e\xi) |u|^2 ) -
\frac{1}{p+1}\int_\Rn K(\e\xi)|u|^{p+1}  .
\end{equation}
Let
\[
\mathcal{M} = \left\{u\in E \colon
\textstyle\int_\Rn (|\n u|^2 + |u|^2 + V(\e\xi)|u|^2 ) 
=\textstyle \int_\Rn |u|^{p+1}  \right\}
\]
be the Nehary manifold of $J$, which has codimension one. Let
\[
\mathcal{N}=\left\{u\in E \colon \textstyle\int_\Rn \left( \left| \left( \frac{\n}{i} - A(\e\xi)
\right) u\right|^2 + |u|^2 + V(\e\xi)|u|^2 \right) 
= \textstyle\int_\Rn |u|^{p+1} \ \right\}
\]
be the Nehari manifold of $F^{\e\xi,\sigma}$.
One checks readily that $\operatorname{codim} \mathcal{N}=1$. 
Recall (\cite{ku}) that $U_c$ is, up to multiplication by a constant phase,
the {\em unique} minimum of $J$ restricted to $\mathcal{M}$. Now, for every
$u\in\mathcal{M}$, the function $x \mapsto e^{i A(\e\xi)\cdot x} u(x)$ lies in
$\mathcal{N}$, and viceversa. Moreover
\[
J(u)=F^{\e\xi,\sigma} (e^{i A(\e\xi)\cdot x} u).
\]
This immediately implies that $\min_{\mathcal{N}} F^{\e\xi,\sigma}$ is
achieved at a point which differs from $e^{i A(\e\xi)\cdot x} U_c (x)$ at most
for a constant phase. In other words, $z^{\e\xi,\sigma}$ is a critical point
for $F^{\e\xi,\sigma}$ of mountain-pass type, and the claim follows by standard
results (see \cite{chang}).

\medskip

Let $R\gg 1$ and consider a radial smooth function 
$\chi_{1}:\R^{n}\longrightarrow \R$ such that 
\begin{equation}\label{eq:c1}  \chi_{1}(x) = 1, \quad \hbox{ for
} |x| \leq R; \qquad  \chi_{1}(x) = 0, \quad \hbox{ for } |x| \geq 2 R; 
\end{equation} 
\begin{equation}\label{eq:c2} 
|\n \chi_{1}(x)| \leq \frac{2}{R}, \quad \hbox{ for } R \leq |x| \leq 2 R. 
\end{equation} 
We also set $ \chi_{2}(x)=1-\chi_{1}(x)$. 
Given $\phi$ let us consider the functions 
$$ 
\phi_{i}(x)=\chi_{i}(x-\xi)\phi(x),\quad i=1,2. 
$$ 
A straightforward computation yields: 
$$ 
\irn |\phi|^2 = \irn |\phi_1|^2 + \irn |\phi_2|^2 + 2\re\irn \phi_{1} \, \bar 
\phi_{2}, 
$$ 
$$ 
\irn |\n \phi|^2 = \irn |\n \phi_1|^2 + \irn |\n \phi_2|^2 + 2\re\irn 
\n\phi_{1} \cdot \overline{\n \phi_{2}}, 
$$ 
and hence 
$$ 
\| \phi \|^2 = \| \phi_1 \|^2 + \| \phi_2 \|^2+ 
2 \re\irn\left[ \phi_{1} \, \overline \phi_{2}
+ \n\phi_{1} \cdot \overline{\n \phi_{2}}\right]. 
$$ 
Letting $I$ denote the last integral, 
one immediately finds: 
$$ 
I=\underbrace{\irn \chi_{1}\chi_{2}(\phi^{2}+|\n \phi|^{2}) }_{I_{\phi}} + 
\underbrace{\irn\phi^{2}\n\chi_{1}\cdot \n\chi_{2}  }_{I'}+ 
\underbrace{\irn (\phi_{2}\n\chi_{1}\cdot\overline{\n\phi}+\overline{\phi_{1}}\,\n 
\phi\cdot\n\chi_{2})}_{I''}. 
$$ 
Due to the definition of $\chi$, the two integrals $I'$ and $I''$ 
reduce to integrals from $R$ and $2R$, and thus they are $o_{R}(1)\|\phi\|^{2}$,
where $o_{R}(1)$ is a function which tends to $0$, as $R \to +\infty$.
As a consequence we have that 
\begin{equation}\label{eq:d} 
\| \phi \|^2 = \| \phi_1 \|^2 + \| \phi_2 \|^2 + 2I_\phi +
o_R(1)\| \phi \|^2 \,.
\end{equation} 
After these preliminaries, let us evaluate the three terms in the 
equation below: 
$$ 
(L_{\e,\sigma,\xi}\phi|\phi)= 
\underbrace{(L_{\e,\sigma,\xi}\phi_{1}|\phi_{1})}_{\a_{1}}+ 
\underbrace{(L_{\e,\sigma,\xi}\phi_{2}|\phi_{2})}_{\a_{2}}+ 
2\underbrace{(L_{\e,\sigma,\xi}\phi_{1}|\phi_{2})}_{\a_{3}}. 
$$ 
One has: 
\begin{multline*}
\alpha_{1}=\langle L_{\e,\sigma,\xi}\phi_{1}\mid \phi_{1}
\rangle =D^2
F^{\e\xi,\sigma}({{z^{\e\xi,\sigma}}})(\phi_{1},\phi_{1}) +
\int_\Rn [V(\e x)-V(\e\xi)]|\phi_{1}|^2  \\ -\int_\Rn [K(\e
x)-K(\e\xi)]|\phi_{1}|^{p+1} + \int_\Rn \bigg|
\bigg( \bigg(\frac{\nabla}{i}-A(\e x)\bigg)\phi_{1}\bigg|^2 - \bigg|
\bigg(\frac{\nabla}{i}-A(\e \xi)\bigg)\phi_{1}\bigg|^2 \bigg) \, .
\end{multline*}

In order to use \eqref{eq:17}, we introduce the function 
$\phi^*_{1}=\phi_{1}-\psi$, where 
$\psi$ is the projection of $\phi_{1}$ onto $\cal V$: 
\begin{multline*}
\psi=(\phi_{1}|z^{\e\xi,\sigma})z^{\e\xi,\sigma} \| z^{\e\xi,\sigma} \|^{-2} + 
(\phi_{1}|iz^{\e\xi,\sigma})iz^{\e\xi,\sigma}\|z^{\e\xi,\sigma}\|^{-2}\\
+\sum (\phi_{1}|\partial_{x_{i}}{{z^{\e\xi,\sigma}}}) 
\partial_{x_{i}}{z^{\e\xi,\sigma}}\| \partial_{x_{i}}{{{z^{\e\xi,\sigma}}}} \|^{-2}. 
\end{multline*}
Then we have: 
\begin{equation} 
D^{2}F^{\e\xi}[\phi_{1},\phi_{1}]= 
D^{2}F^{\e\xi}[\phi_{1}^*,\phi_{1}^*]+ 
D^{2}F^{\e\xi}[\psi,\psi]+2\re D^{2}F^{\e\xi}[\phi^*_{1},\psi] \,.
\label{eq:alfa2} 
\end{equation} 
Since ${z^{\e\xi,\sigma}}$ is orthogonal to $\partial_{x_{i}}{{{z^{\e\xi,\sigma}}}}$, 
$i=1,\ldots,n$, then one 
readily 
checks that $\phi^*_{1}\perp {\cal V}$ and hence \eqref{eq:17} 
implies 
\begin{equation} 
D^{2}F^{\e\xi}[\phi^*_{1},\phi^*_{1}]\geq c_{1} 
\|\phi^*_{1}\|^{2}. 
\label{eq:alfa3} 
\end{equation} 
On the other side, since $(\phi|{{z^{\e\xi,\sigma}}})=0$ it follows: 
\begin{eqnarray*} 
        (\phi_{1}|{{z^{\e\xi,\sigma}}}) & = & (\phi|{z^{\e\xi,\sigma}})-(\phi_{2}|{{z^{\e\xi,\sigma}}})= 
-(\phi_{2}|{z^{\e\xi,\sigma}}) \\ 
& = & -\re\irn\phi_{2}{{z^{\e\xi,\sigma}}}-\re\irn \n {{z^{\e\xi,\sigma}}}\cdot \n \phi_{2} \\ 
&=& -\re\irn\chi_{2}(y)z(y)\phi(y+\xi)dy-\re\irn \n z(y)\cdot \n 
\chi_{2}(y)\phi(y+\xi)dy . 
\end{eqnarray*} 
Since $\chi_{2}(x)=0$ for all $|x|<R$, and since $z(x)\to 0$ as 
$|x|=R\to\infty$, we infer 
$(\phi_{1}|{{z^{\e\xi,\sigma}}})=o_{R}(1)\|\phi\|$.  Similarly one shows that 
$(\phi_{1}|\partial_{x}{{{{z^{\e\xi,\sigma}}}}})=o_{R}(1)\|\phi\|$ and it follows that 
\begin{equation}\label{eq:alfa4} 
         \|\psi\|=o_{R}(1)\|\phi\|. 
\end{equation} 
We are now in position to estimate the last two terms in equation (\ref{eq:alfa2}). 
Actually, using Lemma 3.1 we get 

\begin{eqnarray}\label{eq:alfa5}
& & D^{2}F^{\e\xi}[\psi,\psi] \geq  C\|\psi\|^{2}
+ V(\e\xi)\irn\psi^{2} \\ &- & \notag  K(\e\xi)\bigg[ \re (p-1)
\int_\Rn |z^{\e\xi,\sigma}|^{p-3}
\re (z^{\e\xi,\sigma} \bar \psi) z^{\e\xi,\sigma} \bar \psi  
\\ & & \notag 
+\int_\Rn |z^{\e\xi,\sigma}|^{p-1}|\psi|^2 \bigg] =
o_{R}(1)\|\phi\|^{2}. 
\end{eqnarray}
The same arguments readily imply 
\begin{equation}\label{eq:alfa6} 
\re D^{2}F^{\e\xi}[\phi^*_{1},\psi]=o_{R}(1)\|\phi\|^{2}.
\end{equation} 
Putting together (\ref{eq:alfa3}), (\ref{eq:alfa5}) and  (\ref{eq:alfa6}) 
we infer 
\begin{equation} \label{eq:alfa7} 
D^{2}F^{\e\xi}[\phi_{1},\phi_{1}]\geq C\|\phi_{1}\|^{2}+o_{R}(1)\|\phi\|^{2}.
\end{equation} 
Using arguments already carried out before, one has 
\begin{eqnarray*} 
\irn|V(\e x)-V(\e\xi)|\phi_{1}^{2} & \leq 
&\e c_{2}\irn|x-\xi|\chi_1^{2}(x-\xi)\phi^{2}(x) \\ 
& \leq & \e c_{3}\irn |y|\chi_1^{2}(y)\phi^{2}(y+\xi)dy \\ 
& \leq & \e c_{4} R \|\phi\|^{2},
\end{eqnarray*} 
and similarly for the terms containing $K$.
This and (\ref{eq:alfa7}) yield 
\begin{equation} 
        \a_{1}=(L_{\e,\sigma,\xi}\phi_{1}|\phi_{1})\geq c_{5}
        \|\phi_{1}\|^{2}-
        \e c_{4}R \|\phi\|^{2}+o_{R}(1)\|\phi\|^{2}. 
\label{eq:alfa8} 
\end{equation} 
Let us now estimate $\a_{2}$.  One finds 
\begin{equation} 
\a_{2} = \langle L_{\e,\sigma,\xi} \phi_{2} \mid \phi_{2} \rangle \geq 
c_{6} \|\phi_{2}\|^{2}+o_{R}(1)\|\phi\|^{2}.  \label{eq:alfa9} 
\end{equation} 
In a quite similar way one shows that 
\begin{equation} 
\a_{3} \geq c_{7}I_{\phi}+o_{R}(1)\|\phi\|^{2}. 
\label{eq:alfa10} 
\end{equation} 
Finally, (\ref{eq:alfa8}), (\ref{eq:alfa9}), (\ref{eq:alfa10}) and the fact 
that $I_{\phi}\geq 0$, yield 
\begin{eqnarray*} 
(L_{\e,\sigma,\xi}\phi|\phi) & = & \a_{1}+\a_{2}+2\a_{3}\\ 
& \geq & c_{8}\left[\|\phi_{1}\|^{2}+\|\phi_{2}\|^{2} 
+2 I_{\phi}\right]-c_{9} R \e\|\phi\|^{2}+ o_{R}(1)\|\phi\|^{2}. 
\end{eqnarray*} 
Recalling (\ref{eq:d}) we infer that 
$$ 
(L_{\e,\sigma,\xi}\phi|\phi)\geq c_{10}\|\phi\|^{2}-c_{9} R \e\|\phi\|^{2}+ 
o_{R}(1)\|\phi\|^{2}. 
$$ 
Taking $R = \e^{-1/2}$, and choosing $\e$ small,
equation  \eqref{eq:14} follows. This completes the proof. 
\end{proof}

\section{The finite dimensional reduction}
\label{sec:fdr} 
 In this Section we will show that the existence of critical points of $f_{\e}$ 
 can be reduced to the search of critical points of an auxiliary finite 
 dimensional functional.  The proof will be carried out in two 
 subsections dealing, respectively, with  a Liapunov-Schmidt reduction, and 
 with the behaviour of the auxiliary finite dimensional functional. 

\smallskip

\subsection{A Liapunov-Schmidt type reduction}\label{subsec:LS} The main 
result of this section is the following lemma. 
 
\begin{lemma}\label{lem:w} 
For $\e>0$ small, $|\xi|\leq \overline{\xi}$ and $\sigma\in [0,2\pi]$, there
exists a unique  $w=w(\e,\sigma,\xi)\in 
(T_{z^{\e\xi,\sigma}} Z^\e)^{\perp}$ such that 
$\nabla f_\e (z^{\e\xi,\sigma} + w)\in T_{z^{\e\xi,\sigma}} Z^\e$. 
Such a $w(\e,\sigma,\xi)$ 
 is of class $C^{2}$, resp.  $C^{1,p-1}$, with respect to $\xi$, provided that 
 $p\geq 2$, resp.  $1<p<2$. 
Moreover, the functional $\Phi_\e (\sigma,\xi)=f_\e
(z^{\e\xi,\sigma}+w(\e,\sigma,\xi))$ has  the same regularity as $w$ and
satisfies:   $$ 
\n \Phi_\e(\sigma_0,\xi_0)=0\quad \Longleftrightarrow\quad \n 
f_\e\left(z_{\xi_0}+w(\e,\sigma_0,\xi_0)\right)=0. 
$$ 
\end{lemma} 
\begin{proof} 
Let $P=P_{\e\xi,\sigma}$ denote the projection onto $(T_{z^{\e\xi,\sigma}} Z^\e)^\perp$. We 
want 
to find a 
solution $w\in (T_{z^{\e\xi,\sigma}} Z)^{\perp}$ of the equation 
$P\nabla f_\e(z^{\e\xi,\sigma} +w)=0$.  One has that $\n f_\e(z+w)= 
\n f_\e (z)+D^2 f_\e(z)[w]+R(z,w)$ with $\|R(z,w)\|=o(\|w\|)$, uniformly 
with respect to 
$z=z^{\e\xi,\sigma}$, for  $|\xi|\leq \overline{\xi}$. Using the notation introduced 
in the previous section, we are led to the equation: 
$$ 
L_{\e,\sigma,\xi}w + P\n f_\e (z)+PR(z,w)=0. 
$$ 
According to Lemma \ref{lemma:3.1}, this is equivalent to 
$$ 
w = N_{\e,\xi,\sigma}(w), \quad \mbox{where}\quad 
N_{\e,\xi,\sigma}(w)=-L_{\e,\sigma,\xi}^{-1}\left( P\n f_\e (z)+PR(z,w)\right). 
$$ 
{From} Lemma \ref{lem:1} it follows that 
 
\begin{equation}\label{eq:N} 
        \|N_{\e,\xi,\sigma}(w)\|\leq c_1 (\e|\n V(\e\xi)|++\e|\nabla
K(\e\xi)|+\e |J_A (\e\xi)|+\e^2)+ o(\|w\|).          \end{equation} 
Then one readily checks that $N_{\e,\xi,\sigma}$ is a contraction on some ball in 
$(T_{z^{\e\xi,\sigma}} Z^\e)^{\perp}$ 
provided that $\e>0$ is small enough and $|\xi|\leq \overline{\xi}$. 
Then there exists a unique $w$ such that $w=N_{\e,\xi,\sigma}(w)$.  Let us 
point out that we cannot use the Implicit Function Theorem to find 
$w(\e,\xi,\sigma)$, because the map
$(\e,u)\mapsto P\n f_\e (u)$ fails to be
$C^2$.  However, fixed $\e>0$ small, we can apply the Implicit 
Function Theorem to the map $(\xi,\sigma,w)\mapsto P\n f_\e (z^{\e\xi,\sigma} +
w)$.  Then, in particular, the function $w(\e,\xi,\sigma)$ turns out to be of
class  $C^1$ with respect to $\xi$ and $\sigma$.  Finally, it is a standard
argument, see  \cite{ambad, ABC}, to check that the critical points of
$\Phi_\e (\xi,\sigma)=f_\e (z+w)$ give rise to critical points of $f_\e$. 
\end{proof} 
 
\begin{Remark}\label{remark:psi}
Since $f_\e (z^{\e\xi,\sigma})$ is actually independent of $\sigma$, the implicit function $w$ is 
constant with respect to that variable. As a result, there exists a functional $\Psi_\e \colon \Rn 
\to \R$ such that
\[
\Phi_\e (\sigma,\xi)=\Psi_\e (\xi), \qquad \forall\sigma\in [0,2\pi], \quad \forall \xi\in\Rn.
\]
In the sequel, we will omit the dependence of $w$ on $\sigma$, even it is defined over 
$S^1\times\Rn$.
\end{Remark}

\begin{Remark}\label{rem:w} 
From (\ref{eq:N}) it immediately follows that: 
\begin{equation}        \label{eq:w} 
        \|w\|\leq C \left(\e |\n V(\e\xi)|+\e|\n K(\e\xi)| +\e
|J_A (\e\xi)|+\e^2\right),  \end{equation} 
where $C>0$. 
 \end{Remark} 
 
\noindent  The following result can be proved by adapting the same argument as
in \cite{ams}.

\begin{lemma}\label{lem:Dw} 
One has that: 
\begin{equation}        \label{eq:Dw} 
        \|\nabla_\xi w\|\leq c \left(\e |\n 
V(\e\xi)|+\e|\n K(\e\xi)|+\e |J_A(\e\xi)|+O(\e^2)\right)^\gamma,
\end{equation} 
where $\gamma=\min\{1,p-1\}$ and $c > 0$ is some constant.
\end{lemma}

\subsection{The finite dimensional functional}\label{subs:3}

The purpose of this subsection is to give an explicit form to the
finite--dimensional functional $\Phi_\e (\sigma,\xi)=\Psi_\e (\xi)=f_\e
(z^{\e\xi,\sigma}+w(\e,\xi))$.

Recall the precise definition of $z^{\e\xi,\sigma}$ given in \eqref{def:zeta}.
For brevity, we set in the sequel $z=z^{\e\xi,\sigma}$ and $w =
w(\e,\xi)$.

Since $z$ satisfies \eqref{eqlimite}, we easily find the following relations:
\begin{equation}
\int_\Rn \bigg|\bigg( \frac{\n}{i} - A(\e\xi)\bigg)z \bigg|^2 + |z|^2 +
V(\e\xi)|z|^2=\int_\Rn K(\e\xi)|z|^{p+1}
\end{equation}
\begin{multline}
\re\int_\Rn \bigg( \frac{\n}{i} - A(\e\xi)\bigg)z \cdot \overline{\bigg(
\frac{\n}{i} - A(\e\xi)\bigg)w}  +  \re\int_\Rn z\bar w  \\
+\re\int_\Rn V(\e\xi)z \bar w  = \re\int_\Rn K(\e\xi) |z|^{p-1}z\bar w \,. 
\end{multline}

Hence we get
\begin{multline}\label{eq:Phi}
\Phi_\e (\sigma,\xi)=f_\e (z^{\e\xi,\sigma}+w(\e,\sigma,\xi))
=\\ =K(\e \xi) \left( \frac{1}{2}-\frac{1}{p+1}\right) \int_\Rn
|z|^{p+1}
+ \frac{1}{2} \int_\Rn |A(\e\xi)-A(\e x)|^2 z^2\,  
+\\ \re\int_\Rn (A(\e\xi) -A(\e x))z \cdot (A(\e\xi)-A(\e x))\bar w 
 +\e\re\int_\Rn \frac{1}{i} z \bar w \, \operatorname{div}A(\e x)
\\ + \frac{1}{2}\int_\Rn \bigg| \bigg(\frac{\n}{i} - A(\e x)\bigg)w \bigg|^2 
+\re\int_\Rn [V(\e x)-V(\e\xi)]z \bar w
\\ + \frac{1}{2}\int_\Rn [V(\e x)-V(\e\xi)] |w|^2
+ \frac{1}{2} \int_\Rn [V(\e x)-V(\e\xi)] z^2
\\ +\frac{1}{2} V(\e\xi) \int_\Rn |w|^2\,
- \frac{1}{p+1} \re  \int_\Rn K(\e x)(|z+w|^{p+1} -
|z|^{p+1} - (p+1) |z|^{p-1} z \bar w )\, 
\\ + \re K(\e\xi)\int_\Rn |z|^{p-1}z\bar w 
+O(\e^2) \,.
\end{multline}
Here we have used the estimate
\[
\int_\Rn \bigg(\frac{1}{2} K(\e x)-\frac{1}{p+1} K(\e\xi) \bigg) |z|^{p+1} \,
 = \bigg( \frac{1}{2} - \frac{1}{p+1} \bigg) \int_\Rn K(\e\xi) |z|^{p+1} \,
 + O(\e^2),
\]
which follows from the boundedness of $K''$. Since we know that
\begin{eqnarray*}
\alpha (\e\xi)&=&\left( \frac{1+V(\e\xi)}{K(\e\xi)}\right)^{\frac{1}{p-1}}\\
\beta (\e\xi)&=&\left( 1+V(\e\xi)\right)^\frac{1}{2},
\end{eqnarray*}
we get immediately
\begin{equation}
\int_\Rn |z^{\e\xi,\sigma}|^{p+1}  = C_0 \, \Lambda (\e\xi)
[K(\e\xi)]^{-1},
\end{equation}
where we define the auxiliary function
\begin{equation}\label{def:auxiliary}
\Lambda (x)=\frac{(1+V(x))^\theta}{K(x)^{-2/(p-1)}},\qquad \theta =
\frac{p+1}{p-1}-\frac{n}{2},
\end{equation}
and $C_0 = \|U\|_{L^2}$.
Now one can estimate the various terms in \eqref{eq:Phi} by  means of
\eqref{eq:w} and \eqref{eq:Dw}, to prove that
\begin{equation}\label{eq:expansion}
\Phi_\e (\sigma,\xi) = \Psi_\e (\xi)= C_1 \, \Lambda (\e\xi) + O(\e).
\end{equation}
Similarly,
\begin{equation}\label{eq:expansionD}
\nabla \Psi_\e (\xi )=C_1 \n \Lambda (\e\xi) + \e^{1+\gamma} O(1),
\end{equation}
where $C_1 = \left(
\frac{1}{2}-\frac{1}{p+1} \right)C_0$. We omit the details, which can be
deduced without effort from \cite{ams}.

\section{Statement and proof of the main results}

\bigskip

In this section we exploit the finite-dimensional reduction
performed in the previous section to find existence and multiple
solutions of (NLS).
Recalling Lemma~\ref{lem:w}, we have to look for critical points
of $\Phi_\e$ as a function of the variables
$(\sigma,\xi)\in [0,2\pi ]\times \Rn$
(or, equivalently, $(\eta,\xi)\in S^1 \times \Rn$).

In what follows, we use the
following notation: given a set $\Omega \subset \Rn$ and a number $\rho > 0$,
\[
\Omega_\rho \overset{\rm def}{=} \{ x\in\Rn \mid \operatorname{dist}(x,\Omega) < \rho \}.
\]

We start with the following result, which deals with local extrema.

\begin{theorem}\label{pippo}
Suppose that (K1), (V1) and (A1) hold. Assume moreover that there is
a compact set $M \subset \Rn$ over which $\Lambda$ achieves an isolated strict local 
minimum with value $a$. By this we mean that for some $\delta > 0$,
\begin{equation}\label{eq:slm}
b\; \overset{\rm def}{=} \inf_{x\in \de M_\delta} \Lambda (x) > a.
\end{equation}
Then there exists $\e_\delta > 0$ such that $(S_\e)$ has at least $\operatorname{cat} 
(M,M_\delta)$ (orbits of) solutions concentrating near $M_\delta$, for all $0 < \e < \e_\delta$.

Conversely, assume that $K$ is a compact set $\Rn$ over which $\Lambda$
achieves an isolated strict local maximum with value $b$, namely
for some $\gamma > 0$,
\begin{equation}\label{eq:max}
a \; \overset{\rm def}{=} \inf_{x\in \de K_\gamma} \Lambda (x) < b.
\end{equation}
Then there exists $\e_\gamma > 0$ such that $(S_\e)$ has at
least $\operatorname{cat} (K,K_\gamma)$ (orbits of)
solutions concentrating near $K_\gamma$, for all $0 < \e < \e_\gamma$.

\end{theorem}

\begin{proof}
As in the previous theorem, one has $\Phi_\e (\eta,\xi)=\Psi_\e (\xi)$.
Now choose $\bar\xi > 0$ in such a way that $M_\delta \subset \{x\in\Rn \mid \ |x| < \bar\xi\}$. 
Define again $\overline{\Lambda}$ as in the proof of Theorem \ref{th:5.3}. Let
\begin{eqnarray*}
N^\e &=& \{\xi\in \Rn \mid \e\xi\in M\} \\
N_\delta^\e &=& \{ \xi\in \Rn \mid \e\xi\in M_\delta\}\\
\Theta^\e &=& \{ \xi\in \Rn \mid \Psi_\e (\xi) \leq C_1 \tfrac{a+b}{2}\}.
\end{eqnarray*}
From \eqref{eq:expansion} we get some $\e_\delta > 0$ such that
\begin{equation}\label{eq:inclusioni}
N^\e \subset \Theta^\e \subset N_\delta^\e ,
\end{equation}
for all $0 < \e < \e_\delta$. To apply standard category theory, we need to prove that 
$\Theta^\e$ is compact. To this end, as can be
readily checked, it suffices to prove that $\Theta^\e$ cannot touch $\de N_\delta^\e$. But if 
$\e\xi\in\de M$, one has $\overline{\Lambda} (\e\xi) \geq b$ by the very definition of $\delta$, 
and so
\[
\Psi_\e (\xi) \geq C_1 \overline{\Lambda} (\e\xi) + o_\e (1) \geq C_1 b + o_\e (1).
\]
On the other hand, for all $\xi\in \Theta^\e$ one has also $\Psi_\e (\xi) \leq C_1 \frac{a+b}{2}$. 
We can conclude from \eqref{eq:inclusioni} and elementary properties of the 
Lusternik--Schnirel'man category that $\Psi_e$ has at least
\[
\operatorname{cat} (\Theta^\e,\Theta^\e) \geq \operatorname{cat} (N^\e , N_\delta^\e) = 
\operatorname{cat} (N,N_\delta)
\]
critical points in $\Theta^\e$, which correspond to at least $\operatorname{cat}(M,M_\delta)$ 
orbits of solutions to $(S_\e)$.
Now, let $(\eta^*,\xi^*)\in S^1 \times M_\delta$ a critical point of $\Phi_\e$. Hence this point  
$(\eta^*,\xi^*)$ localizes a solution $u_{\e,\eta^*,\xi^*} (x)=z^{\e\xi^*,\eta^*} (x) + 
w(\e,\eta^*,\xi^*)$ of $(S_\e)$.  Recalling the change of variable which allowed us to pass
from (NLS) to $(S_\e)$, we find that
\[
u_{\e,\eta^*,\xi^*} (x) \approx z^{\e\xi^*,\eta^*} (\tfrac{x-\xi^*}{\e}).
\]
solves (NLS). The concentration statement follows from standard arguments (\cite{ABC,ams}).
The Proof of the second part of Theorem \ref{th:5.3} follows with analogous
arguments.
\end{proof}

Theorem 1.1 in the Introduction is an immediate corollary of the previous one when $x_0$ is 
either a nondegenerate local maximum or minimum for $\Lambda$.
We remark that the case in which $\Lambda$ has a maximum cannot be handled
using a direct variational approach and the arguments in \cite{cingolani}
cannot be applied. 

To treat the general case, we need some more work.
In order to  present our main result, we need to introduce
some topological concepts.

\smallskip
Given a set $M \subset \Rn$, the {\em cup long} of $M$ is by definition
\[
\ell (M)=1+\sup \{k\in\mathbb{N} \mid (\exists \alpha_1,\dots,\alpha_n \in \check{H}^{*} 
(M)\setminus \{1\})(\alpha_1 \cup \dots \cup \alpha_k \neq 0)\}.
\]
If no such classes exists, we set $\ell (M)=1$. Here $\check{H}^{*}(M)$ is the Alexander 
cohomology of $M$ with real coefficients, and $\cup$ denotes the cup product. It is well known 
that $\ell (S^{n-1})=\cat (S^{n-1})=2$, and $\ell (T^n)=\cat (T^n) = n+1$, where $T^n$ is the 
standard $n$--dimensional torus. But in general, one has $\ell (M) \leq \cat (M)$.

The following definition dates back to Bott (\cite{bott}).

\begin{definition}
We say that $M$ is non-degenerate for a $C^2$ function
$I\colon \R^N \to \R$ if $M$ consists of Morse theoretically
non-degenerate critical points for the restriction $I_{|M^\bot}$.
\end{definition}

To prove our existence result, we need the next theorem, which is a slightly modified statement 
of
Theorem 6.4 in chapter II of \cite{chang}.

\begin{theorem}\label{th:chang}
Let $I\in C^1 (V)$ and $J\in C^2 (V)$ be two functionals defined on
the Riemannian manifold $V$,
and let $\Sigma \subset V$ be a smooth, compact, non-degenerate manifold of critical
points of $J$. Denote by $\mathcal{U}$ a neighborhood of $\Sigma$.

If $\| I-J\|_{C^1 (U)}$ is small enough, then the functional $I$ has at least $\ell (\Sigma)$ critical 
points contained in $\mathcal{U}$.
\end{theorem} 

We only remark that Theorem \ref{th:chang} can also be proved in the framework of Conley 
theory (\cite{conley}).

We are now ready to prove an existence and multiplicity result for (NLS).

\begin{theorem}\label{th:5.3}
Let (V1), (K1) and (A1) hold. If the auxiliary function $\Lambda$ has
a smooth, compact, non-degenerate manifold of critical points $M$,
then for $\e > 0$ small, the problem $(S_\e)$ has at least $\ell (M)$ (orbits of) solutions 
concentrating near points of $M$.
\end{theorem}
\begin{proof}
By Remark \ref{remark:psi}, we have to find critical points of
$\Psi_\e=\Psi_\e (\xi)$. Since $M$ is compact, we can choose $\bar\xi > 0$ so
that $|x| < \bar\xi$ for all points $x\in M$. From this moment, $\bar\xi$ is
kept fixed. 
form $\{\eta^*\} \times M$ is obviously  a non-degenerate critical manifold
We set now $V=\Rn$, $J=\Lambda$, $\Sigma = M$, and $I(\xi)=\Psi_\e (\eta,\xi /
\e)$. Select $\delta > 0$ so that $M_\delta \subset \{x \colon |x| < \bar\xi
\}$, and no critical points of $\Lambda$ are in $M_\delta$, except fot those
of $M$. Set $\mathcal{U}=M_\delta$. From \eqref{eq:expansion} and
\eqref{eq:expansionD} it follows that $I$ is close to $J$ in $C^1
(\overline{\mathcal{U}})$ when $\e$ is very small. We can apply Theorem
\ref{th:chang} to find at least $\ell (M)$ critical points
$\{\xi_1,\dots,\xi_{\ell (M)}\}$ for $\Psi_\e$, provided $\e$ is small enough.
Hence the orbits $S^1\times\{\xi_1\}, \dots, S^1\times\{\xi_{\ell (M)}\}$
consist of critical points for $\Phi_\e$ which produce solutions of $(S_\e)$.
The concentration statement follows as in \cite{ams}. \end{proof}

\noindent
\begin{Remark}\label{nonl}   
We point out that
Theorem 1.1, Theorem \ref{pippo} 
and Theorem \ref{th:5.3} hold for problems involving
more general nonlinearities $g(x,u)$
satisfying the same assumptions in \cite{grossi} (see also Remark 5.4
in \cite{ams}).
For our approach, we need the uniqueness of the radial solution
$z$ of the corresponding scalar equation
\begin{equation}\label{eta}
- \varDelta u + u + V(\e \xi)u=g(\e\xi,u), \quad u>0, \ \ u \in W^{1,2}(\Rn) \,.
\end{equation}
Let us also remark that in \cite{cingolani} the class
of nonlinearities handled does not require that
equation \eqref{eta} has a unique solution.
\end{Remark}

\section*{Acknowledgements}
The authors would like to thank Prof.~Ambrosetti for
several comments and suggestions.

\end{article}
\end{document}